\documentclass[12pt]{article}
\usepackage{bez123,calc,curves,ebezier,epic,eepic,graphicx,multiply,rotating}
\everymath{\displaystyle}
\usepackage{graphicx}
\usepackage{pst-all}
\usepackage{color}
\usepackage{amssymb}
\usepackage{amsmath}

\usepackage{epsfig}
\usepackage{pst-grad} 
\usepackage{pst-plot} 

\newtheorem{prethm}{{\bf Theorem}}

\newenvironment{thm}{\begin{prethm}{\hspace{-0.5
               em}{\bf .}}}{\end{prethm}}
\newtheorem{prelemma}{{\bf Lemma}}

\newtheorem{preex}{{\bf Example}}

\newtheorem{preprop}{{\bf Proposition}}

\newtheorem{precor}{{\bf Corollary}}

\newtheorem{preremark}{{\bf Remark}}

\newenvironment{remark}{\begin{preremark}{\hspace{-0.5
               em}{\bf.}}}{\end{preremark}}
\newtheorem{preprob}{{\bf Problem}}

\newenvironment{prob}{\begin{preprob}{\hspace{-0.7
               em}{\bf.}}}{\end{preprob}}
\newtheorem{predefin}{{\bf Definition}}

\newenvironment{defin}{\begin{predefin}{\hspace{-0.5
               em}{\bf .}}}{\end{predefin}}
\newtheorem{preconj}{{\bf Conjecture}}

\newtheorem{preprobb}{{\bf Problem}}

\newtheorem{prelem}{{\bf Theorem}}

\newenvironment{proof}{{\bf Proof.}\rm }{\hfill{$\Box$}}

\newtheorem{presolution}{{\bf Solution.}}

\def\newpic#1{}
\def\qed{\ifhmode\unskip\nobreak\fi\quad\ifmmode\Box\else$\Box$\fi}

\title{\Large\bf\noindent Dynamic monopolies in directed graphs: the spread of unilateral influence in social networks}

\author{\large\bf Kaveh Khoshkhah, Hossein Soltani, Manouchehr Zaker\footnote{Corresponding author: mzaker@iasbs.ac.ir}
\vspace{5mm}\\
   Department of Mathematics,\\
    Institute for Advanced Studies in Basic Sciences,\\
  Zanjan 45137-66731, Iran}
    \date{June 2012}
\begin{document}
\maketitle
\begin{abstract}
\noindent Irreversible dynamic monopolies were already defined and widely studied in the literature for undirected graphs. They are arising from formulation of the irreversible spread of influence such as disease, opinion, adaptation of a new product, etc., in social networks, where the influence between any two individuals is assumed to be bilateral or reciprocal. But in many phenomena, the influence in the underlying network is unilateral or one-sided. In order to study the latter models we need to introduce and study the concept of dynamic monopolies in directed graphs. Let $G$ be a directed graph such that the in-degree of any vertex $G$ is at least one. Let also ${\mathcal{\tau}}: V(G)\rightarrow \Bbb{N}$ be an assignment of thresholds to the vertices of $G$. A subset $M$ of vertices of $G$ is called a dynamic monopoly for $(G,\tau)$ if the vertex set of $G$ can be partitioned into $D_0\cup \ldots \cup D_t$ such that
$D_0=M$ and for any $i\geq 1$ and any $v\in D_i$, the number of edges from $D_0\cup \ldots \cup D_{i-1}$ to $v$ is at least $\tau(v)$. One of the most applicable and widely studied threshold assignments in directed graphs is strict majority threshold assignment in which for any vertex $v$, $\tau(v)=\lceil (deg^{in}(v)+1)/2 \rceil$, where $deg^{in}(v)$ stands for the in-degree of $v$. By a strict majority dynamic monopoly of a graph $G$ we mean any dynamic monopoly of $G$ with strict majority threshold assignment for the vertices of $G$. In this paper we first discuss some basic upper and lower bounds for the size of dynamic monopolies with general threshold assignments and then obtain some hardness complexity results concerning the smallest size of dynamic monopolies in directed graphs. Next we show that any directed graph on $n$ vertices and with positive minimum in-degree admits a strict majority dynamic monopoly with $n/2$ vertices. We show that this bound is achieved by a polynomial time algorithm. This upper bound improves greatly the best known result. The final note of the paper deals with the possibility of the improvement of the latter $n/2$ bound.
\end{abstract}




\section{Introduction and motivation}

\noindent The irreversible spread of influence in social networks such as spread of disease, of opinion and etc. are modeling in terms of progressive (or irreversible) dynamic monopolies in combinatorial undirected graphs \cite{DR,Z}. In this formulation the elements of the network are represented by the nodes of a graph $G=(V,E)$ and the links of the network by the edges of $G$. Assume that corresponding to any vertex $v$ of $G$ an integer value denoted by $\tau(v)$ is given. This value is called the threshold of $v$ and the assignment $v\rightarrow \tau(v)$ is called a threshold assignment of $G$. Let a graph $G$ and an assignment of thresholds $\tau$ to its vertices be given. By a $\tau$-dynamic monopoly we mean any subset $D$ of $G$ such that the vertex set of $G$ can be partitioned into subsets $D_0, D_1, \ldots, D_k$ such that $D_0=D$ and for any $i=1, \ldots, k-1$ each vertex $v$ in $D_{i+1}$ has at least $\tau(v)$ neighbors in $D_0\cup \ldots \cup D_i$. Dynamic monopolies were widely studied in the recent years with various types of threshold assignments and for various families of graphs \cite{ABW,BHLN,CL1,CL3,DR,FKRRS,KKT,KSZ,Z}. The usual formulation of dynamic monopolies is in terms of a discrete time dynamic process defined as follows. Consider a discrete time dynamic process on the vertices of $G$, where some vertices of $G$ are considered as active vertices at the beginning of the process i.e. at time zero (activeness is interpreted according to the underlying phenomenon such as disease, opinion and etc.). Denote the set of active vertices at any discrete time $t\geq 0$ by $D_t$. Assume that at the beginning of the process (i.e. at time zero), the vertices of a subset $D\subseteq V(G)$
are active. Hence $D_0=D$. At each discrete time $i$ any un-active vertex $v$ is activated provided that $v$ has at least $\tau(v)$ active neighbors in $D_0 \cup \ldots \cup D_{i-1}$. If at the end of the process all vertices are active then the starting subset $D$, is called {\it dynamic monopoly} or simply {\it dynamo}. Some well-known threshold assignments for the vertices of a graph $G$ are simple and strict majority thresholds to be defined later.

\noindent While formulating the spread of influence by undirected graphs it is assumed that the influence is a mutual property i.e. when a vertex $v$ does influence another vertex $u$ then $u$ too does influence $v$. We notice that in some applications influence is a unilateral or one-sided relationship. For instance a person may have an influential role
to another person but does not effect from the same person. For such models we have to use directed graphs and extend the concept of dynamic monopolies for directed graphs. Throughout this paper we consider simple directed graphs. A directed graph $G=(V,E)$ is simple if it contains no loop and there exists at most one edge between any two vertices $G$. In particular, there exists no directed cycle of length two in $G$. We refer the reader for other concepts concerning directed graphs not defined in this paper to \cite{W}. Although we consider simple directed graphs but some of our theorems are still valid for multiple directed graphs. We make a remark on this point at the concluding remarks of the paper. We present the following formal definition.

\noindent \begin{defin}
Let $G$ be a directed graph such that the in-degree of any vertex $G$ is at least one. Let also ${\mathcal{\tau}}: V(G)\rightarrow \Bbb{N}$ be an assignment of thresholds to the vertices of $G$ such that $\tau(v)\leq deg^{in}(v)$, for any vertex $v$, where $deg^{in}(v)$ stands for the in-degree of $v$. A subset $M$ of vertices of $G$ is called a dynamic monopoly for $(G,\tau)$ if the vertex set of $G$ can be partitioned into $D_0\cup \ldots \cup D_t$ such that
$D_0=M$ and for any $i\geq 1$ and any $v\in D_i$, the number of edges from $D_0\cup \ldots \cup D_{i-1}$ to $v$ is at least $\tau(v)$.\label{def}
\end{defin}

\noindent For any two vertices $u$ and $v$ if there is an edge from $u$ to $v$ then we say $u$ is an in-neighbor of $v$. Let us remark that since in this model any vertex can only be effected by its in-neighbor vertices then it is assumed that all directed graphs in this paper have positive minimum in-degree. We denote the order of $G$ by $|G|$.

\noindent Two special types of threshold assignments are mostly studies in the area of dynamic monopolies both in directed and undirected graphs. Let $G$ be a directed graph by the simple (resp. strict) majority threshold for $G$ we mean the threshold function $\tau$ such that $\tau(v)=deg^{in}(v)/2$ (resp. $\tau(v)=\lceil (deg^{in}(v)+1)/2 \rceil$) for any vertex $v$ of $G$, where $deg^{in}(v)$ stands for the in-degree of $v$. By a strict majority dynamic monopoly for a graph $G$ we mean any dynamic monopoly for $G$ with strict majority threshold assignment. Strict majority dynamic monopolies were widely studied in the literature. First in \cite{CL1}, Chang and Lyuu have obtained the upper bound $23|G|/27$ for the smallest size of strict majority dynamic monopoly in any directed graph $G$. Then the same authors improved this bound to $0.7732|G|$ in \cite{CL2}. Recently this bound improved to $2|G|/3$ in \cite{CL3} and by a very shorter proof in \cite{ABW} by Ackerman et al. We show in Section 2 of this paper that the smallest size of strict majority dynamic monopoly in any directed graph $G$ is at most $|G|/2$.

\noindent {\bf The outline of the paper is as follows}. In the rest of this section we discuss an upper and a lower bound for the size of dynamic monopolies with general thresholds. Then in Section 2, we obtain some hardness results concerning the complexity status of determining the smallest size of dynamic monopolies with strict majority threshold and with constant threshold assignment $\tau(v)=2$. Next in Section 3, we first show that any strongly connected directed graph admits a strict majority dynamic monopoly with at most $\lceil |G|/2 \rceil$ vertices (Theorem \ref{strong1}). Then we reduce the latter bound to $\lfloor |G|/2 \rfloor$ (Theorem \ref{strong2}). In fact to prove this bound we need the proof of the upper bound $\lceil |G|/2 \rceil$. Finally using this result we show that any directed graph $G$ contains a strict majority dynamic monopoly with at most $\lfloor |G|/2 \rfloor$ vertices (Theorem \ref{main}). Such a strict majority dynamic monopoly can be obtained by a polynomial time algorithm (Remark \ref{rem3}). At the last section we first show that the upper bound of Theorem \ref{main} can not be improved to any bound better than $(2/5)|G|$, i.e. to any bound with order of magnitude $(2/5)|G|-o(1)$. We end the paper with mentioning an open question about the smallest size of strict majority dynamic monopolies.

\noindent For directed graphs with general thresholds we have the following interesting result from \cite{ABW}. Recall that the in-degree of any vertex $v$ is denoted by $deg^{in}(v)$.

\begin{thm} (\cite{ABW})
Let $G$ be any directed graph and $\tau:V(G)\rightarrow \Bbb{N}$ any threshold assignment to the vertices of $G$. Then there exists a dynamic monopoly for $(G,\tau)$ with cardinality at most$$\sum_{v\in V(G)} \frac{\tau(v)}{deg^{in}(v)+1}.$$
\end{thm}

\noindent In order to obtain a lower bound we need the following result from \cite{KSZ}. Let $G$ be a graph and $\tau$ be a threshold assignment to the vertices of $G$. Denote the edge density of $G$ by $\epsilon(G)$. Let also $\bar{t}$ and $t_M$ denote the average and maximum threshold of $\tau$, respectively. It was proved in \cite{KSZ} that for any $\tau$-dynamic monopoly $M$ of $G$ we have $|M|\geq |G|(1-\frac{\epsilon(G)}{\bar{t}})(\frac{\bar{t}}{t_M})$.

\noindent By the similar proof we have the following analogous result for directed graphs.

\begin{thm}
Let $G$ be a directed graph and $\tau$ be a threshold assignment to the vertices of $G$.
Let also $\bar{t}$ and $t_M$ denote the average and maximum threshold of $\tau$, respectively. For any $\tau$-dynamic monopoly $M$ of $G$
we have $$|M|\geq |G|(1-\frac{\epsilon(G)}{\bar{t}})(\frac{\bar{t}}{t_M}).$$
\end{thm}

\section{The complexity results}

\noindent In this section we consider the complexity status of the smallest dynamic monopolies in directed graphs with constant threshold and strict majority threshold. The algorithmic aspects of dynamic monopolies in undirected graphs were studied by various authors (e.g. \cite{DR,C,BHLN, NNUW}). In proving our results for directed graphs we need the results of \cite{C} by Chen, where some hardness and inapproximability results were obtained. We only need the NP-hardness results for dynamic monopolies with strict majority and constant threshold assignments, which we present in the following uniformed format.

\begin{thm}(\cite{C})
To determine the smallest size of dynamic monopolies in undirected graphs is NP-hard either for strict majority threshold and for the case where all vertices have constant threshold 2.\label{chen}
\end{thm}

\noindent In the following two theorems we show that the same results hold for directed graphs.

\begin{thm}
It is an NP-complete problem to determine the smallest size of any dynamic monopoly in directed graphs where the threshold of any vertex is 2.
\end{thm}

\begin{figure}
\begin{center}
\scalebox{1} 
{
\begin{pspicture}(0,-1.883711)(9.46291,1.883711)
\psline[linewidth=0.04cm](2.8010156,1.1837109)(0.6210156,1.2037109)
\psline[linewidth=0.04cm](2.8010156,1.163711)(0.6210156,-0.9962891)
\psline[linewidth=0.04cm](0.6210156,1.1837109)(2.7810156,-0.97628903)
\psline[linewidth=0.04cm](0.6210156,-1.0162891)(2.8010156,-1.0162891)
\psdots[dotsize=0.12](2.8010156,1.163711)
\psdots[dotsize=0.12](0.6210156,1.2037109)
\psdots[dotsize=0.12](0.6210156,-1.0162891)
\psdots[dotsize=0.12](2.8010156,-0.9962891)
\psline[linewidth=0.04cm,arrowsize=0.05291667cm 2.0,arrowlength=1.4,arrowinset=0.4]{->}(1.6810156,1.1837109)(1.8010156,1.1837109)
\psline[linewidth=0.04cm,arrowsize=0.05291667cm 2.0,arrowlength=1.4,arrowinset=0.4]{->}(2.2610157,0.62371093)(2.1210155,0.48371094)
\psline[linewidth=0.04cm,arrowsize=0.05291667cm 2.0,arrowlength=1.4,arrowinset=0.4]{->}(1.3210156,0.48371094)(1.2010156,0.62371093)
\psline[linewidth=0.04cm,arrowsize=0.05291667cm 2.0,arrowlength=1.4,arrowinset=0.4]{->}(2.2810156,-0.45628905)(2.2210157,-0.39628905)
\psline[linewidth=0.04cm,arrowsize=0.05291667cm 2.0,arrowlength=1.4,arrowinset=0.4]{->}(1.2210156,-0.41628906)(1.1810156,-0.45628905)
\psline[linewidth=0.04cm,arrowsize=0.05291667cm 2.0,arrowlength=1.4,arrowinset=0.4]{->}(1.5610156,-1.0162891)(1.8210156,-1.0162891)
\psdots[dotsize=0.12](1.7210156,0.08371094)
\usefont{T1}{ptm}{m}{n}
\rput(2.1424706,0.10871094){$a$}
\usefont{T1}{ptm}{m}{n}
\rput(3.1224706,1.268711){$c$}
\usefont{T1}{ptm}{m}{n}
\rput(0.3124707,1.268711){$b$}
\usefont{T1}{ptm}{m}{n}
\rput(0.2724707,-0.9912891){$u$}
\usefont{T1}{ptm}{m}{n}
\rput(3.1124706,-0.9912891){$v$}
\psline[linewidth=0.04cm](6.6010156,-1.7962891)(8.801016,-1.7962891)
\psline[linewidth=0.04cm](6.6210155,1.8037109)(8.801016,1.8037109)
\psline[linewidth=0.04cm](7.7010155,0.7037109)(8.801016,1.783711)
\psline[linewidth=0.04cm](7.6810155,0.7037109)(6.6210155,1.783711)
\psline[linewidth=0.04cm](7.7010155,-0.6762891)(7.7010155,0.7037109)
\psline[linewidth=0.04cm](7.6810155,-0.69628906)(6.6010156,-1.756289)
\psline[linewidth=0.04cm](7.7010155,-0.69628906)(8.801016,-1.756289)
\psdots[dotsize=0.12](6.6010156,1.8037109)
\psdots[dotsize=0.12](8.801016,1.783711)
\psdots[dotsize=0.12](7.7010155,0.68371093)
\psdots[dotsize=0.12](7.7010155,-0.71628904)
\psdots[dotsize=0.12](6.5610156,-1.7762891)
\psdots[dotsize=0.12](8.821015,-1.7762891)
\psline[linewidth=0.04cm,arrowsize=0.05291667cm 2.0,arrowlength=1.4,arrowinset=0.4]{->}(7.5210156,-1.7962891)(7.821016,-1.7962891)
\psline[linewidth=0.04cm,arrowsize=0.05291667cm 2.0,arrowlength=1.4,arrowinset=0.4]{->}(7.1410155,-1.236289)(7.1010156,-1.2762891)
\psline[linewidth=0.04cm,arrowsize=0.05291667cm 2.0,arrowlength=1.4,arrowinset=0.4]{->}(7.7010155,0.14371094)(7.7010155,-0.09628906)
\psline[linewidth=0.04cm,arrowsize=0.05291667cm 2.0,arrowlength=1.4,arrowinset=0.4]{->}(8.201015,1.1837109)(8.321015,1.283711)
\psline[linewidth=0.04cm,arrowsize=0.05291667cm 2.0,arrowlength=1.4,arrowinset=0.4]{->}(7.861016,1.8037109)(7.6410155,1.8037109)
\psline[linewidth=0.04cm,arrowsize=0.05291667cm 2.0,arrowlength=1.4,arrowinset=0.4]{->}(7.1210155,1.283711)(7.2010155,1.1837109)
\psline[linewidth=0.04cm,arrowsize=0.05291667cm 2.0,arrowlength=1.4,arrowinset=0.4]{->}(8.221016,-1.1762891)(8.1610155,-1.116289)
\usefont{T1}{ptm}{m}{n}
\rput(9.152471,-1.711289){$v$}
\usefont{T1}{ptm}{m}{n}
\rput(6.1924706,-1.711289){$u$}
\usefont{T1}{ptm}{m}{n}
\rput(7.392471,-0.57128906){$x$}
\end{pspicture}
}
\caption{The widget $W_{u,v}$ (left) and the widget $W'_{u,v}$ (right)}\label{fig}
\end{center}
\end{figure}
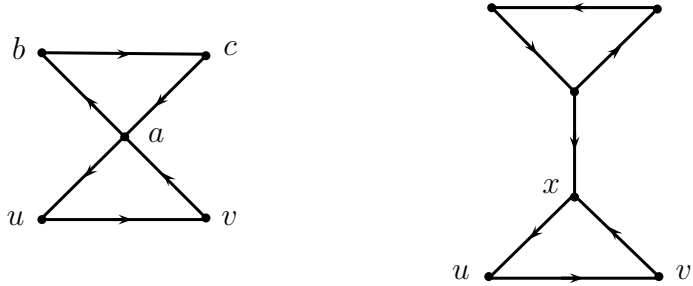

\noindent \begin{proof}
By Theorem \ref{chen} it is enough to obtain a polynomial time reduction from the problem of smallest dynamic monopolies in undirected graphs with constant threshold 2 to our problem. Let $G$ be an undirected graph whose vertices have constant threshold 2. We obtain a directed graph $H$ from $G$ where $deg^{in}(x)\geq 1$ and $\tau(x)=2$ for any $x\in H$. For any edge $uv$ of $G$, we replace the edge $uv$ by the widget $W_{u,v}$, as illustrated in the left side of Figure \ref{fig}. Denote the resulting directed graph by $H$. We may assume that $V(G)\subseteq V(H)$. The graph $H$ consists of $|E(G)|$ widgets $W_{u,v}$ corresponding to any edge $uv$ of $G$. First note that $deg^{in}(x)\geq 1$ for any vertex $x\in H$ and $|V(H)|=|V(G)|+3|E(G)|$. Set $\tau(x)=2$ for any vertex $x$ of $H$. Denote the smallest size of dynamic monopolies in $G$ and $H$ by $d(G)$ and $d(H)$, respectively. In the following we show that $d(H)=d(G)+2|E(G)|$. This implies the assertion of the theorem. Each of such widgets $W_{u,v}$ in $H$ contains a vertex of type $b$ and a vertex of type $c$ as we labeled in the definition of $W_{u,v}$ (see $W_{u,v}$ in Figure \ref{fig}). Note that all of these vertices of type $b$ or $c$ belong to any dynamic monopoly of $H$ since $\tau(b)=\tau(c)=2$ and $deg^{in}(b)=deg^{in}(c)=1$. Let $N$ be the set consisting of all vertices of type $b$ or $c$ in any widget of $H$. Since any dynamic monopoly of $H$ contains $N$, the modified threshold of any vertex of type $a$ is one. But any such vertex has an in-neighbor from the vertices of $G$ (i.e. $v$). It is clear now that if $M$ is any dynamic monopoly for $G$ then $M\cup N$ is a dynamic monopoly for $H$. Conversely, let $K$ be any dynamic monopoly for $H$. We have $N\subseteq K$. If a vertex of type $v$ belongs to $K$ we remove it from $K$ and add the vertex $v$ to $K$. The resulting set $K'$ is still a dynamo for $H$ whose cardinality is not greater than $|K|$. Now, $K'\setminus N$ is dynamo for $G$. This completes the proof.
\end{proof}

\noindent The next theorem deals with the complexity of strict majority dynamic monopolies.

\begin{thm}
It is an NP-complete problem to determine the smallest size of any strict majority dynamic monopoly in directed graphs.
\end{thm}

\noindent \begin{proof}
By Theorem \ref{chen}, it is enough to obtain a polynomial time reduction from the problem of smallest strict majority dynamic monopolies in undirected graphs to our problem. Let $G$ be an undirected graph with strict majority threshold. We obtain a directed graph $H$ from $G$ where $deg^{in}(x)\geq 1$. For any edge $uv$ of $G$, we replace the edge $uv$ by the widget $W'_{u,v}$, as illustrated in Fig. \ref{fig}. Denote the resulting directed graph by $H$. We may assume that $V(G)\subseteq V(H)$. The graph $H$ consists of $|E(G)|$ widgets $W'_{u,v}$ corresponding to any edge $uv$ of $G$. Note that $deg^{in}(x)\geq 1$ for any vertex $x\in H$ and $|V(H)|=|V(G)|+4|E(G)|$. We consider strict majority threshold assignment for the vertices of $H$. Denote the smallest size of strict majority dynamic monopolies in $G$ and $H$ by $dyn(G)$ and $dyn(H)$, respectively. In the following we show that $dyn(H)=dyn(G)+|E(G)|$. This implies the assertion of the theorem. Each widget $W'_{u,v}$ in $H$ consists of two directed triangles with an edge between them; an upper triangle and a lower triangle containing the directed edge $uv$. It is clear that any minimal dynamic monopoly in $H$ contains exactly one vertex from any upper triangle of any widget. After activation of upper triangles the threshold of vertex $x$ in any widget is reduced to one. So $x$ too is activated after activation of the vertex $v$. Also note that $deg_G(u)=deg_H^{in}(u)$ for any vertex $u$ of $G$. It follows that if $M$ is any dynamic monopoly for $G$ then $M\cup N$ is a dynamic monopoly for $H$ where $N$ is any set consisting of only one vertex from each widget $W'_{u,v}$ of $H$. Note that $|M\cup N|=|M|+|E(G)|$.

\noindent Conversely, let $N$ be any dynamic monopoly for $H$. Then $N$ intersects any upper triangle exactly once. If a vertex of type $x$ exists in $N$ we remove it and add the vertex $v$. The resulting set $N'$ is still a dynamo whose cardinality is not greater than $|N|$. Now, $N'\cap V(G)$ is dynamo for $G$. This completes the proof.
\end{proof}

\section{The main results}

\noindent We first list the terminology and notations which are used frequently in the rest of the paper. Assume that a directed graph $G=(V,E)$ is given. In the rest of the paper we denote the vertex and edge set of $G$ by $V$ and $E$, respectively.

\begin{itemize}
\item{For any two vertices $u,v\in V$, by $uv\in E$ we mean there exists an edge from
$u$ to $v$}

\item{By an in-neighbor of a vertex $v$ we mean any vertex $u$ such that $uv\in E$}

\item{$N^{in}(v)=\{u\in V: uv\in E\}$}

\item{The in-degree of $v$ is denoted by $deg^{in}(v)$ and defined as $deg^{in}(v)=|N^{in}(v)|$}

\item{By an ordering $\sigma$ on the vertex set of $G$ on $n$ vertices, we mean any bijective function $\sigma:V\rightarrow \{1, 2, \ldots, n\}$}

\item{Let $\sigma$ be an ordering on the vertex set of $G$. The function $f_{\sigma}:V\rightarrow \Bbb{Z}$ is defined as follows for any vertex $v$ of $G$:
$$f_{\sigma}(v)=|N^{in}(v)\cap \{u: \sigma(u)>\sigma(v)\}|-|N^{in}(v)\cap \{u: \sigma(u)<\sigma(v)\}|$$}

\item{By the $f$-value of any vertex $v$ we mean $f_{\sigma}(v)$}

\item{Let an ordering $\sigma$ and $u,v\in V$ be given. We say $u$ is appeared before $v$ in $\sigma$ (or $v$ is appeared after $u$) whenever $\sigma(u)<\sigma(v)$}

\item{Let an ordering $\sigma$ and $u,v\in V$ be given such that $\sigma(u)<\sigma(v)$. We obtain a new ordering as follows

\begin{equation*}
\sigma'(x) =
\begin{cases}
\sigma(x) & \text{if } \sigma(x)<\sigma(u)~or~\sigma(x)>\sigma(v),\\
\sigma(v) & \text{if } x=u,\\
\sigma(x)-1 & otherwise.
\end{cases}
\end{equation*}

\noindent We say $\sigma'$ is the ordering obtained from $\sigma$ by transmitting the vertex $u$ to after $v$.}

\item{Let an ordering $\sigma$ and $u,v\in V$ be given such that $\sigma(u)<\sigma(v)$. We obtain a new ordering as follows

\begin{equation*}
\sigma'(x) =
\begin{cases}
\sigma(x) & \text{if } \sigma(x)<\sigma(u)~or~\sigma(x)>\sigma(v),\\
\sigma(u) & \text{if } x=v,\\
\sigma(x)+1 & otherwise.
\end{cases}
\end{equation*}

\noindent We say $\sigma'$ is the ordering obtained from $\sigma$ by transmitting the vertex $v$ to before $u$.}

\item{Let $A$ and $B$ be two subsets of $V$. We denote the number of edges from $A$ to $B$ by $d(A,B)$.
If $B=\{v\}$ then we simply write $d(A,v)$.}
\end{itemize}

\noindent In the following we obtain an upper bound for the smallest size of dynamic monopolies in strongly connected directed graphs with
strictly majority thresholds. Recall that a directed graph $G$ is strongly connected if for any two distinct vertices $u$ and $v$ there exists a directed path from $u$ to $v$.

\begin{thm}
\noindent Let $G$ be a strongly connected graph.

\noindent (i) If $G$ contains at least one vertex of in-degree odd, then there exists an ordering $\sigma$ such that for any vertex $v$, $f_{\sigma}(v)\not= 0$.

\noindent (ii) If all in-degrees in $G$ are even then there exists an ordering $\sigma$ such that
$f_{\sigma}(v)\not= 0$ for all but at most one vertex $v$ of $G$.\label{strong}
\end{thm}

\noindent \begin{proof}
We prove the following stronger claim:

\noindent {\bf Claim:} There exists an ordering $\sigma$ satisfying the conditions of the theorem such that, for any $u$, $v$ and $w$ if $f_{\sigma}(u)>0$, $f_{\sigma}(v)<0$ and $f_{\sigma}(w)=0$ then $\sigma(u)< \sigma(w)<
\sigma(v)$.

\noindent We prove the claim by induction on $|G|$. The assertion is trivially hold
when $|G|=1$. Assume that it holds for all graphs of less than $n$ vertices
and let $G$ be a graph with $|G|=n$. If $G$ contains a vertex of in-degree odd then we let $x$ be such a vertex, otherwise let $x$ be an arbitrary vertex of $G$. We remove the vertex $x$ from $G$ and consider the strongly connected components of $G\setminus x$. Let $C$ and $C'$ be
any two components of $G\setminus x$. We note that all the edges between $C$ and $C'$ are either from $C$ to $C'$ or from $C'$ to $C$. We construct a new directed graph $H$ as follows. Corresponding to any strongly connected component $C$ of $G\setminus x$ we consider a vertex say $v_C$.  Let $V(H)=\{v_C :~C~is~a~component~of~G\setminus x\}$. Now for any two components $C$ and $C'$ we put an edge from $v_C$ to $v_{C'}$ if and only if all the edges between $C$ and $C'$ in $G$ are directed from $C$ to $C'$. It is clear that $H$ has no directed cycle. Hence the vertices of $H$ can be ordered as $v_{C_1}, v_{C_2}, \ldots, v_{C_t}$ such that if there exists an edge from $v_{C_i}$ to $v_{C_j}$ in $H$ then $i<j$ (note that this kind of ordering is the same as topological sort of acyclic directed graphs). We order the strongly connected components of $G\setminus x$ according to the very ordering of $V(H)$ i.e. $C_1, C_2, \ldots, C_t$. Note that since $G$ is strongly connected and there exists no edge from $C_2\cup \ldots \cup C_t$ to $C_1$, there exists a vertex say $u\in C_1$ such that $xu\in E$. Now since $|C_1|<|G|$, by the induction hypothesis there exists an ordering for $V(C_1)$ and its corresponding function $f$
such that either no vertex in $C_1$ has $f$-value equal to zero or the only vertex in $C_1$ with zero $f$-value is $u$. Also in $C_1$ the vertices with positive $f$-value are firstly appeared, then the vertex with zero $f$-value and at last the vertices with negative $f$-value are appeared. Now we add the vertex $x$ to the ordering of the vertices of $C_1$ and place it before the negative vertices and after the vertices with positive or zero $f$-value. We obtain an order on $C_1\cup \{x\}$ such that the vertices with zero or positive $f$-value in $C_1$ is now positive (under the new ordering) and the negative vertices in $C_1$ remain negative under the new ordering. Note that at this point we know nothing about the sign of $x$.

\noindent In the previous paragraph we obtained an ordering on $C_1\cup \{x\}$ satisfying the conditions of the theorem. Assume that we have obtained a desired ordering say $\sigma$ on $\{x\}\cup C_1 \cup C_2 \cup \ldots \cup C_{i-1}$, where $i$ is any value with $i\geq 2$. Again, there exists a vertex say $w$ in $C_i$ such that there is an edge from $\{x\}\cup C_1 \cup C_2 \cup \ldots \cup C_{i-1}$ to $w$. By the induction hypothesis $C_i$ admits an ordering such that the only vertex with zero value (if exists)
is the very vertex $w$. Also the vertices with positive value are appeared first and then the vertex with zero value (if exists) and finally the vertices with negative value. Now we place the vertices of $\{x\}\cup C_1 \cup C_2 \cup \ldots \cup C_{i-1}$ before the negative and after the non-negative vertices of $C_i$ such that inside $\{x\}\cup C_1 \cup C_2 \cup \ldots \cup C_{i-1}$ the ordering is the same as $\sigma$. We obtain an order say $\sigma'$ for $\{x\}\cup C_1 \cup C_2 \cup \ldots \cup C_{i-1}\cup C_i$. Since there is no edge from $C_i$ to $\{x\}\cup C_1 \cup C_2 \cup \ldots \cup C_{i-1}$, the positiveness and negativeness of vertices in $\{x\}\cup C_1 \cup C_2 \cup \ldots \cup C_{i-1}$ does not change in $\sigma'$. Therefore all vertices in $\{x\}\cup C_1 \cup C_2 \cup \ldots \cup C_{i-1}\cup C_i$ which are before $x$ take positive $f$-value and all vertices after $x$ take negative $f$-value. For the vertex $x$ we have the following possibilities. If $x$ has odd in-degree then by the definition its $f$-value cannot be zero. If $x$ has an even non-zero in-degree then all other vertices have non-zero $f$-value and all vertices with positive value are before all vertices with negative value. The assertion holds in this case. The only other possibility is when $f_{\sigma}(x)=0$. Since we have taken $x$ arbitrarily then the assertion holds again. This completes the proof.
\end{proof}

\noindent The proof of Theorem \ref{strong} shows the following remark.

\begin{remark}
There exists a polynomial time recursive algorithm which outputs an ordering on the vertex set of $G$ satisfying the conditions of Theorem \ref{strong} for any  strongly connected graph $G$.\label{rem1}
\end{remark}

\noindent Theorem \ref{strong} implies the following bound for strongly connected graphs.

\begin{thm}
Let $G$ be a strongly connected graph on $n$ vertices and set $\tau(v)=\lceil (deg^{in}(v)+1)/2 \rceil$ for any vertex $v$ of $G$. Then the size of smallest dynamic monopoly for $G$ is at most
$\lceil n/2 \rceil$. Moreover, in case that $G$ contains at least one vertex of odd in-degree then $G$ contains a dynamo with no more than $n/2$ vertices.\label{strong1}
\end{thm}

\noindent \begin{proof}
By Theorem \ref{strong} the vertices of $G$ admits an ordering $\sigma$ such that all but at most one vertex has $f$-value non-zero.
Let $M$ be the set of vertices $v$ such that $f_{\sigma}(v)\geq 0$. We observe that $M$ is a strict majority dynamo. In fact the vertices with negative $f$-value become active in turn according to their order in $\sigma$. Similarly if $M$ is the set of vertices $v$ with $f_{\sigma}(v)\leq 0$ then $M$ is a dynamo. The vertices with positive $f$-value become active in turn according to reverse of their order in $\sigma$. Now at least one of these sets has size at most $\lceil n/2 \rceil$. If $G$ contains a vertex of odd in-degree then by Theorem \ref{strong} no vertex has zero $f$-value. Hence
we obtain a dynamo with at most $n/2$ vertices.
\end{proof}

\noindent In the following theorem we improve the bound obtained in Theorem \ref{strong1} for strongly connected graphs.

\begin{thm}
Let $G$ be a strongly connected graph on $n$ vertices and set $\tau(v)=\lceil (deg^{in}(v)+1)/2 \rceil$ for any vertex $v$ of $G$. Then the size of smallest dynamic monopoly for $G$ is at most
$\lfloor n/2 \rfloor$.\label{strong2}
\end{thm}

\noindent \begin{proof}
Assume on the contrary that $G$ does not admit any strict majority dynamic monopoly with $\lfloor n/2 \rfloor$ vertices. Let the ordering $\sigma$ be as in Theorem \ref{strong}. Then by Theorem \ref{strong}, there exists at most one vertex $x$ such that $f_{\sigma}(x)=0$. Set $P=\{v: f_{\sigma}(v)>0\}$ and $N=\{v: f_{\sigma}(v)<0\}$. If $|P|\not= |N|$ then using the construction technique of Theorem \ref{strong1} we obtain a dynamo of size at most $n/2$. If $|P|=|N|$ and $P\cup N=V$ then in this case too the same technique obtains a dynamo with at most $n/2$ vertices. Therefore by Theorem \ref{strong} the only remaining case is when $|V|$ is an odd integer, the in-degree of any vertex is even, $|P|=|N|$ and also there exists a unique vertex $x$ such that $f_{\sigma}(x)=0$. Hence by assumption for $G$ we assume that $G$ and $\sigma$ satisfy the latter properties but $G$ contains no dynamic monopoly with at most $\lfloor n/2 \rfloor$ vertices. We need to define sets $P_1$ and $N_1$ as follows.
$$P_1=\{v\in P:~ d(P\cup \{x\}, \{v\})<d(N, \{v\})\}$$
$$N_1=\{v\in N:~ d(N\cup \{x\}, \{v\})<d(P, \{v\})\}.$$
\noindent Note that $P_1$ in fact consists of the vertices $v$ of $P$ such that if we transmit $v$ to after $x$ in the ordering $\sigma$ then the $f$-value of $v$ corresponding to the new ordering remains positive. The similar property holds for $N_1$. Let also $P_2=P\setminus P_1$ and $N_2=N\setminus N_1$. The arrangement of vertices in $\sigma$ is $P, x, N$ (from left to right). We rearrange the vertices in $P$ (resp. $N$) in such a way that
the vertices of $P_2$ (resp. $N_1$) appear firstly and then come the vertices of $P_1$ (resp. $N_2$) and the relative order of any two vertices in $P_2\cup N_2$ is the same as in the original $\sigma$. Now the order of vertices in the rearranged $\sigma$ is $P_2, P_1, x, N_1, N_2$ (from left to right). We denote this arrangement of $\sigma$ again by $\sigma$.

\noindent In the following we first prove a sequence of facts concerning $P_i$ and $N_i$, $i=1,2$. Then we obtain a contradiction. The contradictions to prove these facts and also the whole theorem are based on the following method. Assume that at some step of the proof, we work with an underlying order $\sigma$ on the vertex set of $G$. If we transform $\sigma$ to another order say $\sigma'$ such that either the cardinality of the related set $P$ corresponding to $\sigma'$ is increased or the cardinality of the related set $N$ corresponding to $\sigma'$ is increased, then similar to the proof of Theorem \ref{strong1}, we obtain a dynamo with no more than $n/2$ vertices, i.e. a contraction.

\noindent {\bf Fact 1.} For any $v\in P_2$ we have $d(P\cup \{x\},v)=d(N,v)$ and for any $v\in N_2$ we have
$d(N\cup \{x\},v)=d(P,v)$.

\noindent {\bf Proof.} By definition $d(P\cup \{x\},v)\geq d(N,v)$ for any $v\in P_2$. Assume on the contrary that there exists $v\in P_2$ such that $d(P\cup \{x\},v)>d(N,v)$. We obtain $\sigma'$ from $\sigma$ by transmitting $v$ to after $x$. Then the size of the related set $N$ corresponding to $\sigma'$ is strictly greater than the size of previous $N$ for $\sigma$. As we mentioned before, this implies the existence of a dynamo with at most $n/2$ vertices, i.e. is a contradiction. For the second part of the lemma, assume that
there exists $v\in N_2$ with $d(N\cup \{x\},v)>d(P,v)$. In this case we consider the ordering obtained from $\sigma$ by transmitting $v$ to before $x$. The rest of the proof is similar.

\noindent An interpretation of Fact 1 is that if we transmit $v$ to after $x$ in $\sigma$ and obtain $\sigma'$ then the new $f$-value of $v$ becomes zero. Similarly, if we transmit $v$ to before $x$ then the new $f$-value of $v$ becomes zero.

\noindent {\bf Fact 2.} $d(P_1, x)=0$, $d(N_1, x)=0$.

\noindent {\bf Proof.} Assume on the contrary that there exists $u\in P_1$ such that there is an edge from $u$ to $x$. Let $\sigma'$ be the ordering obtained from $\sigma$ by transmitting $u$ to after $x$. Then the size of the related $P$ is increased strictly, a contradiction. The proof for $d(N_1, x)=0$ is similar.

\noindent {\bf Fact 3.} $d(P_1, N_2)=0$, $d(N_1, P_2)=0$.

\noindent {\bf Proof.} Assume on the contrary that there exists an edge from some vertex $a\in P_1$ to some vertex $b\in N_2$. We obtain $\sigma'$ from $\sigma$ as follows. First we transmit $b$ to before $x$ and then transmit $a$ to after $b$. We obtain $P_2, P_1\setminus \{a\}, b, a, x, N_1, N_2\setminus \{b\}$ (the vertices are ordered from left to right). For the resulting order by Fact 1 the corresponding set $P$ is increased strictly with respect to the previous $P$. This is a contradiction. The proof for $d(N_1, P_2)=0$ is similar.

\noindent {\bf Fact 4.} For any vertex $u\in N_1$, there exists $v\in P_1$ such that $uv \in E$ and also
$f_{\sigma}(v)=2$.

\noindent {\bf Proof.} Assume on the contrary that such a vertex does not exist. We obtain a new ordering $\sigma'$ from $\sigma$ by transmitting $u$ to the beginning of the ordering. By the definition, the $f$-value of $u$ becomes positive. By our hypothesis the sign of $P_1$ remains unchanged. Also by Fact 3, the $f$-value of any vertex of $P_2$ remains unchanged. Hence the cardinality of the corresponding set $P$ (related to $\sigma'$) is strictly increasing with respect to the previous $P$. This contradiction proves Fact 4.

\noindent {\bf Fact 5.} For any vertex $u\in P_1$, there exists $v\in N_1$ such that $uv \in E$ and also
$f_{\sigma}(v)=-2$.

\noindent {\bf Proof.} The proof is similar to the proof of Fact 4.

\noindent {\bf Fact 6.}
Let $u_1$ and $u_2$ be any two vertices from $N_1$. Then there does not exist $v\in P_1$ such that $u_1v,u_2v\in E$ and $f_{\sigma}(v)=2$.

\noindent {\bf Proof.} Assume on the contrary that such a vertex $v\in P_1$ exists. In $\sigma$ we transmit $v$ to after $x$ and $u_1$ and $u_2$ to before $x$ and obtain $\sigma'$. In fact the arrangement of vertices in $\sigma'$ is as follows (where the vertices are ordered from left to right)
$$P_2, P_1\setminus \{v\}, u_1, u_2, x, v, N_1\setminus \{u_1,u_2\}, N_2$$
\noindent In $\sigma'$, the sign of $f$-value of any vertex of $N$ is the same as before. The sign of $f_{\sigma'}(v)$ is negative. Therefore the size of corresponding $N$ is strictly increased, a contradiction.

\noindent {\bf Fact 7.}
Let $u_1$ and $u_2$ be any two vertices from $P_1$. Then there does not exist $v\in N_1$ such that $u_1v\in E$, $u_2v\in E$ and $f_{\sigma}(v)=-2$.

\noindent {\bf Proof.} The proof is similar to the proof of Fact 6.

\noindent {\bf Fact 8.}
There exists a directed matching from $N_1$ to $P_1$ which saturates $N_1$ and also a directed matching from $P_1$ to $N_1$ which saturates $P_1$.

\noindent {\bf Proof.} The existence of directed matching from $N_1$ to $P_1$ is obtained by Fact 4 and Fact 6. Similarly, Fact 5 and Fact 7 imply the existence of directed matching from $P_1$ to $N_1$.

\noindent The following is obtained from Fact 8.

\noindent {\bf Fact 9.}
There exists a directed matching from $N_1$ to $P_1$ which saturates $N_1\cup P_1$ and also a directed matching from $P_1$ to $N_1$ which saturates $P_1\cup N_1$.

\noindent We note that there is an edge from $P$ to $x$, since otherwise $f(x)>0$. By Fact 2 there is no edge from $P_1$ to $x$. Hence there is an edge from $P_2$ to $x$. Assume that $y\in P_2$ is such that $yx\in E$
and for any $y'\in P$ with $y'x\in E$, one has $\sigma(y)>\sigma(y')$. We claim that there exists a $z\in N_1\cup N_2$ such that $yz\in E$, $f_{\sigma}(z)=-2$ and in $\sigma$, $z$ is the first vertex having these properties. If such a vertex does not exist then by transmitting $y$ to the end of $\sigma$ we obtain $\sigma'$ with $f_{\sigma'}(y)<0$ and since there exists no vertex $z$ with the above-mentioned properties, then for each $u\in N_1\cup N_2$, $f_{\sigma'}(u)<0$. It follows that the number of vertices with negative $f$-value is increased, a contradiction. Therefore the claim is proved and such a vertex $z$ exists. For $z$ there are two possibilities.

\noindent {\bf Case 1.} $z\in N_1$.

\noindent By Fact 9, there exists $u\in P_1$ with $uz\in E$. We obtain $\sigma'$ as follows.

\noindent $P_2\setminus \{y\}~~~P_1\setminus \{u\}~~~z~u~x~y~~~N_1\setminus \{z\}~~~N_2$

\noindent We already had $f_{\sigma}(z)=-2$. We have now edges from $y$ and $u$ to $z$. Hence in $\sigma'$ we have $f_{\sigma'}(z)>0$. We already had $f_{\sigma}(x)=0$ also $yx\in E$ and by Fact 2, there is no edge from $z$ to $x$. Hence in $\sigma'$ we have $f_{\sigma'}(x)>0$. By the definition of $P$, for any $w\in P_1\cup P_2\setminus \{u,y\}$ we have $f_{\sigma'}(w)>0$. Also since $uz\in E$, there is no edge from $z$ to $u$ in $G$. By the definition of $P_1$, $f_{\sigma'}(u)>0$. In other words the number of vertices with positive $f$ is strictly increased. This contradiction completes the proof in this case.

\noindent {\bf Case 2.} $z\in N_2$ and $zx\not\in E$.

\noindent In this case we construct $\sigma'$ as follows.

\noindent $P_2\setminus \{y\}~~~P_1~~~z~x~y~~~N_1~~~N_2\setminus \{z\}$

\noindent By Fact 1 and that $yz\in E$ we have $f_{\sigma'}(z)=+2$.
Since $yx\in E$ and $zx\not\in E$ and also $f_{\sigma}(x)=0$ we have $f_{\sigma'}(x)=+2$.
This implies that the number of vertices with positive $f$-value is increased in $\sigma'$, a contradiction.

\noindent {\bf Case 3.} $z\in N_2$ and $zx\in E$.

\noindent In this case we choose $w\in P_1\cup P_2$ satisfying the following properties.

(i) $f(w)=+2$

(ii) $zw\in E$

(iii) $w$ is the last vertex in $\sigma$ satisfying (i) and (ii).

\noindent If there does not exist such a vertex $w$ we obtain the ordering $\sigma'$ using $\sigma$ by transmitting the vertex $z$ to the beginning of $\sigma$. And then we get contradiction. If $w\in P_1$ we continue similar to case 1. If $w\in P_2$ then there are two possibilities (note that $w$ could not be $y$ because there is no multiple edge).

\noindent Subcase 1. $\sigma(y)<\sigma(w)$.

\noindent In this subcase, since $y$ is the closest vertex to $x$ in $\sigma$ with $yx\in E$ then $wx\not\in E$. The rest of the argument is similar to Case 2.

\noindent Subcase 2. $\sigma(w)<\sigma(y)$.

\noindent In this subcase we construct $\sigma'$ using $\sigma$ by transmitting the vertex $z$ to before $y$. In this new ordering, the $f$-value of $z$ becomes positive by Fact 1 and $yz\in E$. But we had already $f_{\sigma}(z)<0$. The sign of vertices in $P$ remains unchanged (i.e. positive). Now the number of vertices
with positive $f$-value is increased in $\sigma'$.
This contradiction completes the proof.
\end{proof}

\noindent By Remark \ref{rem1}, an ordering $\sigma$ satisfying the conditions of Theorem \ref{strong} can be obtained in a polynomial time. We note that any of the facts in the proof of Theorem \ref{strong2} can be easily checked by a polynomial time procedure. We have therefore the following remark.

\begin{remark}
There exists a polynomial time algorithm which obtains a strict majority dynamic monopoly with no more that $|G|/2$ vertices in any strongly connected graph $G$.\label{rem2}
\end{remark}

\noindent The next theorem deals with general directed graphs. We use the fact that any directed graph can be easily decomposed into vertex disjoint strongly connected components. See \cite{W} for basics in theory of directed graphs.

\begin{thm}
Let $G$ be a directed graph such that no vertex of $G$ has in-degree zero. Let $C_1, C_2, \ldots, C_k$ be the strongly connected components of $G$. For any vertex $v\in C_i$ define $\tau_i(v)=\lceil (deg^{in}_{C_i}(v)+1)/2 \rceil$. Let also $D_i$ be a dynamo for $(C_i, \tau_i)$. Then $D_1 \cup D_2 \cup \ldots \cup D_k$ is a dynamo for $(G, \tau)$ where $\tau(v)=\lceil (deg^{in}_{G}(v)+1)/2 \rceil$. In case that $|C_i|=1$ for some $i$, then we may choose $D_i$ as the null set. \label{pre-last}
\end{thm}

\noindent \begin{proof}
Corresponding to $G$ and $C_1, \ldots, C_k$ we construct a directed graph $H$ as follows.
Corresponding to any component $C_i$ we consider a vertex denoted by $v_i$. Set $V(H)=\{v_1, v_2, \ldots, v_k\}$. We put an edge from $v_i$ to $v_j$ if and only if there exists an edge from $C_i$ to $C_j$ in $G$. It is clear that $H$ is acyclic. We may assume that $v_1, \ldots, v_k$ are arranged such that $v_jv_i\not\in H$ whenever $i<j$. We consider the corresponding order of $v_1, \ldots, v_k$ for $C_1, \ldots, C_k$ too. Note that the in-degree of any vertex $v$ in $C_1$ is the same as the in-degree of $v$ in $G$. This implies that $D_1$ activates all vertices of $C_1$ in $G$. Assume that $D_1\cup \ldots \cup D_{i-1}$ have activated the vertices of $C_1\cup \ldots \cup C_{i-1}$. In the following we show that $C_1\cup \ldots \cup C_{i}$ is activated by $D_1\cup \ldots \cup D_{i}$. To prove this fact we need only to show that the vertices of $C_i$ too are activated. Let $v$ be any vertex of $C_i$. If there exists no edge from $V\setminus C_i$ to $v$ then $\tau_i(v)=\tau(v)$. Otherwise the edges from $V\setminus C_i$ to $v$ are necessarily from $C_1, \ldots, C_{i-1}$ to $v$. Assume that there are say, $r$ edges from $C_1, \ldots, C_{i-1}$ to $v$. In this case $\tau(v)\leq \tau_i(v)+\lceil r/2 \rceil$. But the point is that after activation of all vertices in $C_1\cup \ldots \cup C_{i-1}$ including those in-neighbors of $v$ in $C_1\cup \ldots \cup C_{i-1}$, the threshold of $v$ practically reduces to a value no more than $\tau_i(v)$. Therefore the vertices of $C_i$ are activated according to the activation process in $(C_i, \tau)$ using the dynamo $D_i$. Hence the whole graph $G$ is activated by $D_1 \cup D_2 \cup \ldots \cup D_k$.

\noindent To complete the proof let us note that when for some $i$, $C_i$ consists of only one vertex, then there is no need to put $C_i$ in the dynamic monopoly for $G$. Since in this case all in-neighbors of $C_i$ are in $C_1\cup \ldots \cup C_{i-1}$. After activation of $C_1\cup \ldots \cup C_{i-1}$, $C_i$ is also activated. This completes the proof.
\end{proof}

\noindent Theorem \ref{strong2} and Theorem \ref{pre-last} imply the following result.

\begin{thm}
Let $G$ be a directed graph on $n$ vertices with no vertex of in-degree zero. Set $t(v)=\lceil (deg^{in}(v)+1)/2 \rceil$ for any vertex $v$ of $G$. Then the smallest size of dynamic monopolies for $G$ is at most $\lfloor n/2 \rfloor$.\label{main}
\end{thm}

\noindent As we mentioned in Remark \ref{rem2}, in any strongly connected graph $H$ a strict majority dynamic monopoly of size $|H|/2$ is obtained by a polynomial time algorithm. Also it is easy to decompose any directed graph into its strongly connected components. Theorem \ref{pre-last} and \ref{main} imply the following.

\begin{remark}
There exists a polynomial time algorithm which outputs a strict majority dynamic monopoly with cardinality at most $|G|/2$ for any directed graph $G$ with positive minimum in-degree.\label{rem3}
\end{remark}

\section{Concluding remarks}

\noindent In this section we first show that the upper bound of Theorem \ref{main} can not be improved to any bound better than $(2/5)|G|$, i.e. to any bound with order of magnitude $(2/5)|G|-o(1)$. Next we present an open question
concerning an upper bound for strict majority dynamic monopolies in directed graphs.

\noindent By a $2$-regular directed complete graph on $5$ vertices we mean an edge orientation of the complete graph on $5$ vertices $K_5$, where the in-degree of any vertex is $2$. Since the edge set of $K_5$ is decomposed into two edge disjoint cycles of length $5$, it is easy to obtain a $2$-regular directed $K_5$. In fact orient the edges of either of $5$-cycles and transform each of them into a directed $5$-cycle. The result is a $2$-regular directed $K_5$. Now consider $k$ vertex disjoint copies of $2$-regular $K_5$. Add an extra vertex say $x$ and exactly one edge from each $K_5$ to $x$ so that $deg^{in}(x)=k$. Denote the resulting graph by $G_k$ and the smallest size of any strict majority dynamic monopoly of $G_k$ by $dyn(G_k)$. It is obvious that any dynamo of $G_k$ needs at least two vertices from each copy of $K_5$. Therefore $dyn(G_k)\geq 2k$. In fact equality holds in the latter inequality. We have now
$$\lim_{k\rightarrow \infty} \frac{dyn(G_k)}{|G_k|}=\frac{2}{5}.$$

\noindent Based on the above fact we pose the following question.
\begin{prob}
Is it true that any directed graph $G$ on $n$ vertices with positive minimum in-degree, contains a strict majority dynamic monopoly
with $\lfloor 2n/5 \rfloor$ vertices?
\end{prob}

\noindent We note by Theorem \ref{main} that if the above question is affirmative for strongly connected graphs then it is also affirmative for all directed graphs. We end the paper with the following remark. As we mentioned in Introduction, we considered simple directed graphs in this paper. We note that Theorem \ref{strong1} still holds for multiple directed graphs. But the result of Theorem \ref{main} is not valid for multiple directed graphs. For example consider complete directed graph $\overrightarrow{K_n}$ on $n$ vertices, where $n$ is odd and between any two vertices $u$ and $v$ there exists an edge from $u$ to $v$ and another edge from $v$ to $u$. Observe that any strict majority dynamic monopoly for this graph needs at least $(n+1)/2$ vertices.


\end{document}